\documentclass[11pt]{article}

\usepackage{mathptmx}

\usepackage{amsmath}

\usepackage{amssymb}

\usepackage{amsthm}

\usepackage[all,tips]{xy}

\usepackage{hyperref}

\theoremstyle{plain}

\newtheorem{Step}{Step}

\newtheorem{thm}{Theorem}[section]

\newtheorem*{thmnl}{Theorem}

\newtheorem{prop}[thm]{Proposition}

\newtheorem{lemma}[thm]{Lemma}

\newtheorem{defn}[thm]{Definition}

\theoremstyle{definition}

\newtheorem*{exampleA}{Example A}

\newtheorem*{exampleB}{Example B}

\newtheorem{ques}[thm]{Question}

\newtheorem*{rem}{Remark}

\DeclareMathOperator{\Fix}{Fix}

\DeclareMathOperator{\SL}{SL}

\DeclareMathOperator{\PSL}{PSL}

\DeclareMathOperator{\Ort}{O}

\DeclareMathOperator{\PO}{PO}

\DeclareMathOperator{\PSO}{PSO}

\DeclareMathOperator{\PU}{PU}

\DeclareMathOperator{\SU}{SU}

\DeclareMathOperator*{\vol}{vol}

\DeclareMathOperator{\tr}{Tr}

\DeclareMathOperator*{\Isom}{Isom}

\newcommand{\eps}{\varepsilon}

\newcommand{\abs}[1]{\left\vert#1\right\vert}

\newcommand{\set}[1]{\left\{#1\right\}}

\newcommand{\innp}[1]{\left< #1 \right>}

\newcommand{\pr}[1]{\left( #1 \right) }

\newcommand{\lra}{\longrightarrow}

\newcommand{\B}[1]{\ensuremath{\mathbf{#1}}}

\newcommand{\BB}[1]{\ensuremath{\mathbb{#1}}}

\newcommand{\Cal}[1]{\ensuremath{\mathcal{#1}}}

\newcommand{\Fr}[1]{\ensuremath{\mathfrak{#1}}}

\newcommand{\Hy}{{\mathbf H}}

\newcommand{\Q}{\ensuremath{\BB{Q}}}

\newcommand{\R}{\ensuremath{\BB{R}}}

\newcommand{\Z}{\ensuremath{\BB{Z}}}

\newcommand{\C}{\ensuremath{\BB{C}}}

\begin{document}



\title{\textbf{Length and eigenvalue equivalence}}

\author{C. J. Leininger\thanks{Partially supported by the N. S. F.},
D. B. McReynolds\thanks{Partially supported by a C.M.I. lift-off},
W. D. Neumann\thanks{Partially supported by the N. S. F.} and
A. W. Reid\thanks{Partially supported by the N. S. F.}}

\maketitle



\begin{abstract}
Two Riemannian manifolds are called eigenvalue equivalent when their
   sets of eigenvalues of the Laplace-Beltrami operator are equal
   (ignoring multiplicities). They are (primitive) length equivalent
   when the sets of lengths of their (primitive) closed geodesics are
   equal. We give a general construction of eigenvalue equivalent and
   primitive length equivalent Riemannian manifolds. For example we
   show that every finite volume hyperbolic $n$--manifold has pairs of
   eigenvalue equivalent finite covers of arbitrarily large volume
   ratio. We also show the analogous result for primitive length
   equivalence.
\end{abstract}



\section{Introduction}

Let $M$ be a compact Riemannian manifold, and let $\Delta =
\Delta_M$ denote the Laplace--Beltrami operator of $M$ acting on
$\textrm{L}^2(M)$.  The \textbf{eigenvalue spectrum} ${\cal E}(M)$
consists of the eigenvalues of $\Delta$ listed with their
multiplicities. Two manifolds $M_1$ and $M_2$ are said to be
\textbf{isospectral} if ${\cal E}(M_1)={\cal E}(M_2)$.  Geometric
and topological constraints are forced on isospectral manifolds; for
example if the manifolds are hyperbolic (complete with all sectional
curvature equal to $-1$) then they must have the same volume
\cite{Pe}, and so for surfaces the same genus.

Another invariant of $M$ is the \textbf{length spectrum}
${\cal L}(M)$ of $M$; that is the set of all lengths of closed
geodesics on $M$ counted with multiplicities.  Two manifolds $M_1$
and $M_2$ are said to be \textbf{iso-length spectral} if ${\cal
L}(M_1) = {\cal L}(M_2)$. Under the hypothesis of negative sectional
curvature the invariants ${\cal E}(M)$ and ${\cal L}(M)$ are closely
related. For example, it is known that ${\cal E}(M)$ determines the
set of lengths of closed geodesics, and in the case of closed
hyperbolic surfaces, the stronger statement that ${\cal E}(M)$
determines ${\cal L}(M)$ and vice-versa holds \cite{Hu,Hu2}.

In this paper we address the issue of how much information
is lost by forgetting multiplicities. More precisely, for a compact
Riemannian manifold $M$, define the \textbf{eigenvalue set} (resp.\
\textbf{length set} and \textbf{primitive length set}) to be the set
of eigenvalues of $\Delta$ (resp.\ set of lengths all closed geodesics
and lengths of all primitive closed geodesics) counted without
multiplicities. These sets will be denoted $E(M)$, $L(M)$ and $L_p(M)$
respectively.  Two manifolds $M_1$ and $M_2$ are said to be
\textbf{eigenvalue equivalent} (resp.\ \textbf{length equivalent} and
\textbf{primitive length equivalent}) if $E(M_1) = E(M_2)$ (resp.\
$L(M_1)=L(M_2)$ and $L_p(M_1) = L_p(M_2)$). Although length spectrum
and primitive length spectrum determine each other, the corresponding
statement for length sets is false.  Primitive
length equivalent manifolds are clearly length equivalent, but we shall
see that the converse is false.

We will focus mainly on hyperbolic manifolds of finite volume. Even in
this setting little seems known about the existence of manifolds which
are eigenvalue (resp.\ length or primitive length) equivalent but not
isospectral or iso-length spectral. Examples of non-compact arithmetic
hyperbolic 2--manifolds that are length equivalent were constructed in
Theorem 2 of \cite{Sc} using arithmetic methods.
However, as far as the authors are aware, no examples of closed
hyperbolic surfaces that are length equivalent and not iso-length
spectral were known, and it would appear that no examples of
eigenvalue equivalent or primitive length equivalent hyperbolic
manifolds which are not isospectral or iso-length spectral were known.
Our main results rectify this situation for hyperbolic surfaces and
indeed for all
finite volume hyperbolic $m$--manifolds.

\begin{thm}\label{T:GenLength}
Let $M$ be a closed hyperbolic $m$--manifold.  Then there
exist infinitely many pairs of finite covers $\{M_j,N_j\}$ of $M$
such that
\begin{itemize}
\item[(a)] $E(M_j) = E(N_j)$,
\item[(b)] $\vol(M_j)/\vol(N_j) \to \infty$.
\end{itemize}
Moreover, $E(M_j) = E(N_j)$
for any Riemannian metric on $M$.
\end{thm}

The method of proof of Theorem \ref{T:GenLength} does not
provide (primitive) length equivalent pairs of covers.  However, we
can prove an analogue for primitive length equivalence (and hence also
length equivalence).

\begin{thm}\label{T:GenPrim}
Let $M$ be a finite volume hyperbolic $m$--manifold.  Then there
exist infinitely many pairs of finite covers $\{M_j,N_j\}$ of $M$
such that
\begin{itemize}
\item[(a)] $L_p(M_j) = L_p(N_j)$
\item[(b)] $\vol(M_j)/\vol(N_j) \to \infty$.
\end{itemize}
Moreover, $L_p(M_j) = L_p(N_j)$ for any Riemannian metric on $M$.
\end{thm}

Indeed, as we point out in \S~\ref{se:5.1}, for every finite volume
hyperbolic $n$--manifold where $n\neq 3,4,5$ we can produce pairs
of finite sheeted
covers of arbitrarily large volume ratio that are
both primitive length equivalent and eigenvalue
equivalent.

The methods of the paper are largely group theoretic,
relying on the fundamental group rather than the geometry, and a
quick way to provide lots of examples in many more situations is
given by the following. Recall that a group $\Gamma$ is called
\textbf{large} if it contains a finite index subgroup that surjects a
free non-abelian group.

\begin{thm}\label{T:1} (Theorem \ref{T:Gen2} is a stronger version.)
Let $M$ be a compact Riemannian manifold with large fundamental group.
Then there exist infinitely many pairs
of finite covers $\{M_j,N_j\}$ of $M$ such that
\begin{itemize}
\item[(a)] $L(M_j) = L(N_j)$,
\item[(b)] $E(M_j) = E(N_j)$,
\item[(c)] $\vol(M_j)/\vol(N_j) \to \infty$.
\end{itemize}
Moreover, (a) and (b) hold for any Riemannian metric on $M$, and if
$\pi_1(M)$ is hyperbolic,
$L_p(M_j)=L_p(N_j)$ also holds for any Riemannian metric on $M$.
\end{thm}

Our arguments start with Sunada's construction \cite{Sn} of
isospectral manifolds, which was based on a well known construction in
number theory of ``arithmetically equivalent'' number fields (see
\cite{Per}). Our length equivalence of manifolds similarly has a
number theoretic counterpart called ``Kronecker equivalence'' of
number fields, as we discovered after doing this work; see the
book \cite{Kl}. The results contained here can thus be viewed as
providing the geometric investigation proposed in the sentence from
the last paragraph of that book: \emph{``In view of the relations
between arithmetical and Kronecker equivalence, one should also
study Kronecker equivalence in this geometric situation.''}

In the final section we collate some remarks and questions. In
particular, we note that Mark Kac's famous paper ``Can one hear
the shape of a drum'' \cite{kac} has been a catalyst for much of the work
on isospectrality, and we revisit that paper and the
Gordon--Webb--Wolpert answer to his question \cite{drums} in
the light of our work.



\section{Equivalence} \label{S:Sunada}

We first recall Sunada's construction.

For any finite group $G$ and subgroups $H$ and $K$ of $G$,
we say that $H$ and $K$ are \textbf{almost conjugate} (or ``Gassmann equivalent'' in the
terminology of Perlis \cite{Per}) if for any $g$ in $G$ the
following condition holds (where $(g)$ denotes conjugacy class):

\begin{equation*}
\abs{H \cap (g)} = \abs{K \cap(g)}.
\end{equation*}

In \cite{Sn} Sunada proved the following theorem relating almost
conjugate pairs with isospectral covers.

\begin{thmnl}[Sunada]
Let $M$ be a closed Riemannian manifold, $G$ a finite group, and $H$
and $K$ almost conjugate subgroups of $G$. If $\pi_1(M)$ admits a
homomorphism onto $G$, then the finite covers $M_H$ and $M_K$
associated to the pullback subgroups of $H$ and $K$ are isospectral.
Moreover, the manifolds $M_H$ and $M_K$ are iso-length spectral.
\end{thmnl}

The proof of this is an easy exercise, but checking when
manifolds produced by Sunada's method are non-isometric requires more
work. However, for length equivalence far less is required, and the
resolution of the isometry problem is built into our construction.

Length equivalence, primitive length equivalence, and
eigenvalue equivalence each require a different condition on
the group $G$.  In each instance, we
describe a group theoretic condition, and then explain
how it is used to produce examples with the desired
features.

\subsection{Length and primitive length equivalence}

Though it is not essential, the group $G$ will always
be finite in what follows.

\begin{defn}[Elementwise conjugacy]\label{CC}
Subgroups $H$ and $K$ of $G$ are said to be \textbf{elementwise
conjugate} if for any $g$ in $G$ the following condition holds:

\begin{equation}\label{club}
H \cap (g) \neq \emptyset\quad\hbox{if and only if}\quad K \cap (g) \neq\emptyset.
\end{equation}

(Or, more briefly,$H^G=K^G$, where $H^G=\bigcup_{g\in G}g^{-1}Hg$.)\end{defn}

It is immediate from the definition
that almost conjugate subgroups are elementwise conjugate.

To produce primitive length equivalent manifolds, we impose further
conditions on  $H$ and $K$, and also on $\pi_1(M)$.

\begin{defn}[Primitive]\label{PCC}
We shall call a subgroup $H$ of $G$ \textbf{primitive in $G$} if the
following holds:
\begin{itemize}
\item[(a)] All non-trivial cyclic subgroups of $H$ have the same order $p$ (necessarily prime).
\item[(b)] $\bigcap_{g\in G}g^{-1}Hg=\set1$.
\end{itemize}
\end{defn}

\begin{thm}\label{T:Swallow}
Let $M$ be a Riemannian manifold, $G$ a group, and $H$ and $K$
elementwise conjugate subgroups of $G$.
\begin{itemize}
\item[(1)] If $\pi_1(M)$ admits a homomorphism onto $G$,
then the covers $M_H$ and $M_K$ associated to the pullback subgroups of $H$ and
$K$ are length equivalent.
\item[(2)] If, in addition, $H$ and $K$ are primitive in $G$ and
$\pi_1(M)$ has the property that any pair of distinct maximal cyclic
subgroups of $\Gamma$ intersect trivially, then the covers $M_H$ and
$M_K$ associated to the pullback subgroups of $H$ and $K$ are
primitive length equivalent.
\end{itemize}
\end{thm}

\textbf{Remark.} It is well known that when $M$ admits a
metric of negative sectional curvature, then $\pi_1(M)$ satisfies the
condition needed to apply Theorem \ref{T:Swallow}.

\begin{proof}[Proof of theorem] To prove (1) it suffices to show
that a closed geodesic $\gamma$ on $M$ has a lift to a closed geodesic
on $M_H$ if and only if it has a lift to a closed geodesic on $M_K$.
Let $\rho$ denote the homomorphism $\pi_1(M)\to G$. By standard
covering space theory, $\gamma$ has a closed lift to $M_H$ if and only
if $\rho([\gamma]) \in G$ is conjugate into $H$. By assumption this is
true for $H$ if and only if it is true for $K$, proving (1).

For (2) we will show the inclusion $L_p(M_H) \subseteq L_p(M_K)$;
the reverse inclusion then follows by symmetry. We argue by
contradiction, assuming there is a primitive $\gamma$ in $\pi_1(M_H)$,
every conjugate of which in $\pi_1(M_K)$ is imprimitive.  Let $\gamma_K$ be
any conjugate of $\gamma$ in $\pi_1(M_K)$ and let $\delta \in
\pi_1(M_K)$ and $r > 1$ be such that $\delta^r = \gamma_K$.  The
arguments splits into two cases.









\textbf{Case 1.} $\rho(\delta)=1$. Since $\ker \rho < \pi_1(M_H)$, all conjugates of $\delta$ are
contained in $\pi_1(M_H)$. This contradicts the primitivity of
$\gamma$, as a $\pi_1(M)$--conjugate of $\delta$ powers to $\gamma$.

\textbf{Case 2.}  $\rho(\delta)\ne 1$, so $\rho(\delta)$ has prime
order $p$ by Definition \ref{PCC} (a).  We split this into two subcases.

\textbf{Case 2.1.} $\rho(\gamma_K) \ne 1$.
Since $\innp{\rho(\gamma_K)}$ is nontrivial and contained in
$\innp{\rho(\delta)}$ which has prime order, it is equal to
$\innp{\rho(\delta)}$. Thus, $\rho(\mu\delta\mu^{-1}) \in
\innp{\rho(\gamma)}$, where $\mu$ is the element of $\pi_1(M)$
conjugating $\gamma_K$ to $\gamma$.  Since the cyclic subgroup
$\innp{\rho(\gamma)}$ is contained in $H$, $\rho(\mu\delta\mu^{-1})$
is contained in $H$, so $\mu\delta\mu^{-1}$ is an element of
$\pi_1(M_H)$. This contradicts the primitivity of $\gamma$ since
$(\mu\delta\mu^{-1})^r=\gamma$.

\textbf{Case 2.2.} $\rho(\gamma_K) = 1$.
By Definition \ref{PCC} (b), there exists an element $g$ in $G$ which
conjugates $\rho(\delta)$ outside of $K$.  For any element $\sigma\in
\rho^{-1}(g)$, $\sigma \gamma_K\sigma^{-1}\in \pi_1(M_K)$ and by
assumption this cannot be primitive. Therefore, there exists
$\delta_1$ in $\pi_1(M_K)$ and $s>1$ such that $\sigma \gamma_K
\sigma^{-1} = \delta_1^s$.
We have the equality
$(\sigma\delta\sigma^{-1})^r = \delta_1^s$. By assumption,
$\innp{\sigma\delta\sigma^{-1}}$ and $\innp{\delta_1}$ are contained
in a common maximal cyclic subgroup $C$ of $\pi_1(M)$. The
intersection of $\rho(C)$ with $K$ is a cyclic subgroup which contains
the image of $\innp{\rho(\delta_1)}$. By Definition
\ref{PCC} (a), the cyclic subgroups of $K$ have prime order $p$, and
so $\abs{\rho(C) \cap K} = 1~\hbox{or}~ p$.

Assume first that the latter holds. Now $\rho(C)$ has a unique cyclic
subgroup of order $p$, so $\rho(C)\cap K$ must equal
$\innp{\rho(\sigma\delta\sigma^{-1})}$. Hence the element
$\rho(\sigma\delta\sigma^{-1})$ is in $K$, which contradicts the
choice of $\sigma$. Therefore, we can assume that $\rho(C) \cap K =
1$. Then $\rho(\delta_1)=1$.  Replacing $\gamma_K$ by
$\sigma\gamma_K\sigma^{-1}$ and $\delta$ by $\delta_1$, Case 1
provides the desired contradiction.\end{proof}



\subsection{Eigenvalue equivalence}

To give context to our construction below of eigenvalue equivalence,
we begin by recalling the following well known equivalent formulation of
almost conjugacy.

\begin{prop}\label{SunadaRep}
Subgroups $H$ and $K$ of a finite group $G$ are almost
conjugate if and only if for every finite dimensional complex 
representation $\rho$ of G,
\[ \dim \Fix(\rho|H) = \dim \Fix(\rho|K)\,, \]
where $\Fix(\rho|H)$ denotes the subspace of $\rho(H)$--fixed vectors.
\end{prop}

\begin{proof} A convenient reference for the character theory used here
and later is \cite{Se2}. The dimension of $\Fix(\rho|H)$ is the
inner product $(\chi^H_1,\chi^H_{\rho|H})$ of the trivial character
on $H$ and the character of $\rho|H$. By definition this is
$\frac{1}{\abs{H}}\sum_{h \in H} \chi^H_{\rho|H}(h)$, and since
characters are constant on conjugacy classes, this equals
$\frac1{\abs H}\sum|(g)\cap H|\,\chi^G_{\rho}(g)$, where the sum is
over conjugacy classes in $G$.  Thus the equality $\dim
\Fix(\rho|H)=\dim \Fix(\rho|K)$ is equivalent to
\begin{equation}\label{eq:cc}
\frac1{|H|}\sum|(g)\cap H|\chi^G_{\rho}(g)=\frac1{|K|}\sum|(g)\cap K|\chi^G_{\rho}(g)\,.
\end{equation}
Clearly, almost conjugacy of $H$ and $K$ implies this equality.

For the converse, note first that the equality $\dim
\Fix(\rho|H)=\dim \Fix(\rho|K)$ applied to the regular
representation of $G$ becomes $[G:H]=[G:K]$, whence $|H|=|K|$.  Since
characters of irreducible representations form a basis for class
functions on $G$, letting $\rho$ run over all irreducible
representations of $G$ in equation \eqref{eq:cc} now implies that
$|(g)\cap H|=|(g)\cap K|$ for each conjugacy class $(g)$.
\end{proof}


\begin{defn} We say subgroups $H$ and $K$ of a finite group $G$ are
\textbf{fixed point equivalent} if for any finite dimensional complex
representation $\rho$ of $G$, the restriction $\rho|_H$ has a
nontrivial fixed vector if and only if $\rho|_K$ does.
\end{defn}

\begin{thm}\label{T:UnitSwallow}
Let $H$ and $K$ be fixed point equivalent subgroups of a finite
group $G$.  If $M$ is a compact Riemannian manifold and $\pi_1(M)$
admits a homomorphism onto $G$, then the covers $M_H$ and $M_K$
associated to the pullbacks in $\pi_1(M)$ of $H$ and $K$ are
eigenvalue equivalent.
\end{thm}

\begin{proof}
Let $\widetilde{M}$ be the cover of $M$ associated to the pullback
in $\pi_1(M)$ of the trivial subgroup of $G$.  The action of $G$ on
$\widetilde{M}$ is by isometries, and the quotients by  $H$ and $K$
give covering maps $p_H\colon\widetilde{M} \to M_H$ and
$p_K\colon\widetilde{M} \to M_K$, respectively.

The covering projection induces an embedding $p_H^*\colon
\textrm{L}^2(M_H) \to \textrm{L}^2(\widetilde{M})$ whose image is
the $H$--fixed subspace $\textrm{L}^2(\widetilde{M})^H$ of
$\textrm{L}^2(\widetilde M)$.  Since $\Delta_{\widetilde{M}} \circ
p_H^* = p_H^* \circ \Delta_{M_H}$, the action of $H$ on
$\textrm{L}^2(\widetilde M)$ restricts to an action on the
$\lambda$--eigenspace $\textrm{L}^2(\widetilde M)_\lambda$ and
$p_H^*$ identifies $\textrm{L}^2(M_H)_\lambda$ with
$(\textrm{L}^2(\widetilde M)_\lambda)^H$. Thus $\lambda$ is an
eigenvalue for $M_H$ if and only if $(\textrm{L}^2(\widetilde
M)_\lambda)^H$ has positive dimension.

Since $G$ is finite, the representation of $G$ on
$\textrm{L}^2(\widetilde M)_\lambda$ decomposes as a direct sum of
finite dimensional representations (in fact, compactness of $M$
implies $\textrm{L}^2(\widetilde M)_\lambda$ is finite dimensional,
but we do not need this). Hence, if $H$ and $K$ are fixed point
equivalent, $(\textrm{L}^2(\widetilde M)_\lambda)^H$ will be
non-trivial if and only if $(\textrm{L}^2(\widetilde M)_\lambda)^K$ is
non-trivial.
\end{proof}

\begin{rem} 1. The compactness assumption on $M$ is not necessary. If
$M$ is non-compact our
argument extends easily to show that under the conditions of
the theorem both the discrete and
non-discrete spectra of $M_H$ and $M_K$ agree when viewed as sets.

2. What makes the Sunada construction work for both the length and
eigenvalue spectra is the equivalence of almost conjugacy with the
condition of Proposition \ref{SunadaRep}. Our weakening of almost
conjugacy to elementwise conjugacy on the one hand, and, via
Proposition \ref{SunadaRep}, to fixed point equivalence on the
other, go in dual directions. They therefore cannot be expected to
be equivalent, and it is a little surprising that in the examples we
know, the two weaker conditions still tend to have significant
overlap.
\end{rem}



\subsection{Examples}\label{se:2.3}

An elementary example of elementwise conjugacy is the following.

Let $G$ be the alternating group $\textrm{Alt}(4)$, and
set $a = (12)(34)$ and $b = (14)(23)$. Then the subgroups $H =
\set{1, a}$ and the Klein 4--group $K = \set{1, a, b, ab}$ are
elementwise conjugate.

However note that $K$ is not primitive since it is
a normal subgroup. In addition $H$ and $K$ are not fixed point
equivalent since $K$ has no fixed vector under the irreducible
$3$--dimensional representation of $G$ while $H$ has a fixed
vector.

On the other hand, it is not hard to check that $H$ is fixed point
equivalent to the trivial subgroup $\set1$.

We now generalize this example.

Let $\BB{F}_ p$ be the prime field of order $p$, and let $n\geq 2$ be
a positive integer. The \textbf{$n$--dimensional special
$\BB{F}_p$--affine group} is the semidirect product $\BB{F}_ p^n
\rtimes \SL(n;\BB{F}_p)$ defined by the standard action of
$\SL(n;\BB{F}_p)$ on $\BB{F}_ p^n$. We call any
$\BB{F}_p$--vector subspace $V$ of $\BB{F}_ p^n$ a
\textbf{translational subgroup} of $\BB{F}_p^n \rtimes
\SL(n;\BB{F}_p)$.

\begin{thm}\label{L:Swallow}
Let $V$ and $W$ be translational subgroups of  $G=\BB{F}_p^n \rtimes
\SL(n;\BB{F}_p)$. Then,
\begin{itemize}
\item[(i)] if\/ $V$ and $W$ are both non-trivial then they are
elementwise conjugate in $G$, and they are moreover primitive if they
are proper subgroups of $\BB{F}_p^n$;
\item[(ii)] if\/ $V$ and $W$ are both proper subgroups of $\BB{F}_p^n$
then they are fixed point equivalent in $G$.
\end{itemize}
\end{thm}

\begin{proof} (i). Since $\SL(n;\BB{F}_p)$ acts transitively on
non-trivial elements of $\BB{F}_p^n$, the elementwise conjugacy is
immediate. Moreover, conditions (a) and (b) of Definition \ref{PCC} clearly hold for $V$ if $V$ is a proper
subgroup of $\BB{F}_p^n$.

(ii). It suffices to show that any proper translational subgroup $V$ is
fixed point equivalent to the trivial subgroup. So we must show that
for any $m$-dimensional representation $\rho$ of $\BB{F}_p^n \rtimes
\SL(n;\BB{F}_p)$ with $m>0$, the restriction $\rho|_V$ has a
nontrivial fixed subspace when restricted to $V$.  To this end, let
$\chi$ be the character of $\rho$. The dimension of the fixed space of
$\rho|_V$ is 
$
\dim(\Fix(\rho|_V)) = \frac{1}{\abs{V}} \sum_{v \in V} \chi(v).
$ 
Since $\chi(1)=m$ and any two nontrivial elements of $V$ are
conjugate in $\BB{F}_p^n \rtimes \SL(n;\BB{F}_p)$,
we can rewrite this:
\[ \dim(\Fix(\rho|_V)) = \frac{1}{\abs{V}} \pr{m + (\abs{V} - 1)\chi(x)},\]
where $x \in V-\set0$.
Similarly, the dimension of the
fixed space
for the full translation subgroup $\BB{F}_p^n$ is
\[ \dim(\Fix(\rho|_{\BB{F}_p^n}))=
\frac{1}{\abs{\BB{F}_p^n}} \pr{m + (\abs{\BB{F}_p^n}
  - 1)\chi(x)}\,. \]
Thus $m + (\abs{\BB{F}_p^n} - 1)\chi(x)\ge 0$, so $\chi(x)\ge \frac{-m}{\abs{\BB{F}_p^n}
- 1}$. Hence, $\dim(\Fix(\rho|_V)) = \frac1{\abs V}\left(m + (\abs{V} -
1)\chi(x)\right)\ge \frac1{\abs V}\left(m-m\frac{\abs V-1}{\abs
{\BB{F}_p^n}- 1}\right)>0$.
\end{proof}



\section{Proofs of main results}

The following is a
stronger version of Theorem \ref{T:1}:

\begin{thm}\label{T:Gen2}
Let $M$ be a compact Riemannian manifold whose fundamental group is large. For
every integer $n\geq 2$ and every odd prime $p$,
there exists a finite tower of covers of $M$
\[ M_0 \lra M_{1} \lra \dots \lra  M_{n-1} \lra M_n\lra M\,, \]
with each $M_i\to M_{i+1}$  of degree $p$, such that:
\begin{itemize}
\item[(a)] $L(M_j) = L(M_k)$\quad for $0\le j,k\le n-1$;
\item[(b)] $E(M_j) = E(M_k)$\quad for $1\le j,k\le n$;
\end{itemize}
Moreover, (a),(b) hold for any Riemannian metric on $M$. Finally, if
$\pi_1(M)$ is hyperbolic then  for any Riemannian metric on $M$,
\begin{itemize}
\item[(c)]  $L_p(M_j)=L_p(M_k)$\quad for $1\le j,k\le n-1$.
\end{itemize}
\end{thm}

\begin{proof} Since $M$ is large we can find finite index subgroups which
surject any finitely generated free group, so there is a finite cover
$X$ of $M$ with
\[ \xymatrix{ \pi_1(X) \ar@{->>}[r] & \BB{F}_p^n \rtimes \SL(n;\BB{F}_p)}. \]
Consider any complete $\BB{F}_p$--flag
\[\set0 =V_0\subset V_1 \subset V_2 \subset \dots \subset V_{n-1} \subset V_n = \BB{F}_p^n \]
in $\BB{F}_p^n$. Pulling these subgroups back to $\pi_1(X)\subset
\pi_1(M)$ we obtain a tower
\[ M_0 \lra M_{1} \lra \dots \lra M_{n-1} \lra M_n \lra M\] of
corresponding covers of $M$. The theorem then follows from Theorem
\ref{L:Swallow} combined with Theorems \ref{T:Swallow} and
\ref{T:UnitSwallow}.
\end{proof}

This theorem implies Theorems \ref{T:GenLength}, \ref{T:GenPrim}, and
\ref{T:1} in the case of closed hyperbolic surfaces. In addition, it
is well-known that closed and finite volume hyperbolic manifolds whose
fundamental groups are large exist in all dimensions (see e.g.,
\cite{Lub2}).  This provides examples of hyperbolic manifolds in all
dimensions satisfying the conclusions of Theorems \ref{T:GenLength},
\ref{T:GenPrim}, and \ref{T:1}.  To prove that {\bf any} closed or
finite volume hyperbolic manifold has finite sheeted covers with
these properties requires additional work.

We mention in passing that Theorem \ref{T:Gen2} (b)
applied to surfaces produces arbitrarily long towers of abelian
covers
\[ M_{1} \lra \dots  \lra M_{n-1}\lra M_n \]
whose first nontrivial eigenvalue remains constant. On the other hand,
it is well known that any infinite tower of abelian covers of
a fixed hyperbolic surface has $\lambda_1$ tending to zero (see
\cite{brooks} and \cite{Sn}).


\subsection{More families}

\begin{thm}\label{otherpair}
\begin{itemize}
\item[(i)] Let $p$ be an odd prime.
Then $\PSL(2;{\Z}/p^2{\Z})$
contains subgroups $K<H$ with $[H:K]=p$ which are fixed point equivalent.
\item[(ii)] Let $k\ne\Q$ be a number field with ring of integers
$\Cal{O}_k$.  Let $\Cal P$ be the set of non-dyadic prime ideals
$\Fr{p}$ of $\Cal{O}_k$\, for which $\Cal{O}_k/\Fr p=\BB F_q$ is a
non-prime field (this set is infinite by the Cebotarev Density
Theorem).  Then for $\Fr{p}$ in $\Cal P$ the group
$\PSL(2;\Cal{O}_k/\Fr{p}^2)$ contains subgroups $K<H$ with $[H:K]=p$
which are primitive and elementwise conjugate in
$\PSL(2;\Cal{O}_k/\Fr{p}^2)$.
\end{itemize}
\end{thm}

The proof of Theorem \ref{otherpair} will
be deferred until \S~\ref{se:4}.  Assuming this we complete the proofs of
Theorems \ref{T:GenLength} and \ref{T:GenPrim} in the next subsection.



\subsection{Completion of proofs}

We shall need the following special case of the Strong Approximation
Theorem (see \cite{We} and \cite{No}; see also \cite{LR} for a
discussion of the proof in the particular case of hyperbolic
manifolds). Suppose $M^m$ is a finite volume hyperbolic manifold with
$m\ge 3$.  We shall identify
$\Isom(\Hy^m)$ with $\PO_0(m,1)$ so $\pi_1(M)<\PO_0(m,1)$.  We can
assume there is a number field $k$ such that $\pi_1(M) < \PO_0(m,1;S)$
for a finite extension ring $S$ of $\Cal{O}_k$ with $k$ the field of
fractions of $S$ (see \cite{Rag} for the details). We choose $k$
minimal.

\begin{thm}[Strong Approximation]\label{th:strong approx}
For all but finitely many primes $\frak{p}$ of $S$ the image of
$\pi_1(M)$ under the reduction homomorphism
\[r_{\Fr{p}^j}\colon  \PO_0(m,1;S) \lra \PO(m,1;S/\frak{p}^j)\] 
contains the commutator subgroup $\Omega(m,1;S/\frak{p}^j)$ of $\PO(m,1;S/\frak{p}^j)$ for
all $j\ge 1$.\qed
\end{thm}

\begin{proof}[Proof of Theorem \ref{T:GenLength}]
Theorem \ref{T:GenLength} is shown for hyperbolic surfaces in the
comment following the proof of Theorem \ref{T:Gen2}, so we can
assume that $M$ is a closed hyperbolic manifold of dimension $m\geq
3$.  We will produce surjections of $\pi_1(M)$
onto finite groups containing $\PSL(2;\Z/p^2\Z)$ for infinitely
many $p$.

Let $S$ be as in the Strong Approximation Theorem above.  For all
but a finite number of primes $\frak p$ of $S$ the image of
$\pi_1(M)<\PO_0(m,1;S)$ under the restriction homomorphism
$r_{\Fr{p}^j}$ contains $\Omega(m,1;S/\frak{p}^j)$ for all $j \geq
1$.  If $p$ is the integer prime that $\frak p$ divides then
$r_{\Fr{p}^j}(\pi_1(M))$ therefore contains the subgroup
$\Omega(m,1;\Z/p^j\Z)$, and therefore also the subgroup
$\Omega(2,1;\Z/p^j\Z)$ of $\Omega(m,1;\Z/p^j\Z)$.

We claim the finite groups $\Omega(2,1;\Z/p^j\Z)$ and
$\PSL(2;\Z/p^j\Z)$ are isomorphic.  To see this, first recall that
the $p$--adic Lie groups $\PSL(2;\Q_{\! p})$ and $\Omega(2,1;\Q_{\!
p})$ are isomorphic (as $p$--adic Lie groups).  The groups
$\PSL(2;\Z_{\! p})$ and $\Omega(2,1;\Z_{\! p})$ are the respective
maximal compact subgroups of $\PSL(2;\Q_{\! p})$ and
$\Omega(2,1;\Q_{\! p})$, and are unique up to isomorphism (see
\cite[Ch 3.4]{PR}).  Hence the groups $\PSL(2;\Z_{\! p})$ and
$\Omega(2,1;\Z_{\! p})$ are isomorphic as $p$--adic Lie groups.
Reducing modulo the ideal generated by the $j$th power of the
uniformizer $\pi$ of $\Z_{\! p}$ yields the asserted isomorphism
between $\Omega(2,1;\Z/p^j\Z)$ and $\PSL(2;\Z/p^j\Z)$.

Restricting now to $j=2$ we have shown $r_{\Fr{p}^2}(\pi_1(M))$
contains a subgroup isomorphic to $\PSL(2;\Z/p^2\Z)$.  So by passage
to a subgroup of finite index in $\pi_1(M)$, we can arrange a finite
cover $X$ of $M$ with a surjection $\pi_1(X) \to \PSL(2;\Z/p^2\Z)$.
The existence of pairs $\{M_j,N_j\}$ as stated in Theorem
\ref{T:GenLength} now follows from Theorem \ref{T:UnitSwallow},
Theorem \ref{otherpair}(i), and the infinitude of $\Cal{P}$.
\end{proof}

\begin{proof}[Proof of Theorem \ref{T:GenPrim}]
We have already shown this for hyperbolic surfaces in the comments
following Theorem \ref{T:Gen2}. We next consider hyperbolic
$3$--manifolds.  Let $M$ be a hyperbolic $3$--manifold with holonomy
representation $\pi_1(M) < \PSL(2;S)$ (again $S$ chosen minimally).
The field of fractions of $S$, a finite ring extension of $\Z$, is
necessarily a proper extension of $\Q$, see \cite{MR}.  In
particular, by the Cebotarev Density Theorem, there exist infinitely
many prime ideals $\Fr{p}$ of $S$ such that $S/\Fr{p}$ is a
nontrivial extension of $\BB{F}_p$. The Strong Approximation Theorem
applies here to see that for all but finitely many among this
infinite set of prime ideals of $S$ the reduction maps
\[ r_{\Fr{p}^2}\colon \pi_1(M) \lra \PSL(2;S/\Fr{p}^2), \]
are onto. By Theorem \ref{otherpair} (ii) and Theorem
\ref{T:Swallow} (ii), there exists a pair of covers
\[ N_{\Fr{p}} \lra M_{\Fr{p}} \lra M \] with $L_p(N_{\Fr{p}}) =
L_p(M_{\Fr{p}})$ and $\vol(N_{\Fr{p}})/\vol(M_\Fr{p})=p$.

We extend this to all hyperbolic $m$--manifolds with $m>3$ as follows.
Let $S=\Z[i]$ and let $\Cal P$ be the set of prime ideals
defined in Theorem \ref{otherpair} (specifically, these are the ideals
$p\Z[i]$ with $p\equiv3$ mod 4). For $m>3$ and $\frak p \in
\Cal P$ we first claim we have an injection of
$\PSL(2;\Z[i]/\Fr{p}^j)$ into $\Omega(m,1;\Z/p^j\Z)$.  For this, we argue
as follows. First, there exists a quadratic form $B_4$ defined over
$\Q$ of signature $(3,1)$ and an injection
\[ \PSL(2;\Z[i]) \lra \PSO_0(B_4;\Z). \]
For each prime $\Fr{p} = p\Z[i]$, this induces isomorphisms
\[ \PSL(2;\Z[i]/\Fr{p}^j) \lra \Omega(B_4;\Z/p^j\Z). \]
For $j=1$, this can be found in \cite{Su}. For $j>1$, this is
established by an argument similar
to that used in the proof of the equivalence of $\PSL(2;\Z/p^j\Z)$
and $\Omega(2,1;\Z/p^j\Z)$ in proving Theorem \ref{T:GenLength}.
Extending the form $B_4$ from $\Q^4$ to $\Q^{m+1}$ for $m>3$ by the
identity produces injections
\[ \Omega(B_4;\Z/p^j\Z) \lra \Omega(m,1;\Z/p^j\Z). \]
In particular, we can view
\[ \PSL(2;\Z[i]/\Fr{p}^j) < \Omega(m,1;\Z/p^j\Z) \] for all $m>3$, all
$j$, and all $\frak p\in \Cal P$ as claimed. Since we have
already shown in the proof of Theorem \ref{T:GenLength} that
$\pi_1(M)$ surjects finite groups containing $\Omega(m,1;\Z/p^j\Z)$
for all but finitely many primes, Theorem \ref{otherpair} (ii) with
$S=\Z[i]$ and Theorem \ref{T:Swallow} (ii) now complete the proof.\end{proof}

\section{Proof of Theorem \ref{otherpair}}\label{se:4}


Throughout this section $p$ will be an odd prime. For any ring $R$
let $M(2;R)$ be the algebra of $2\times 2$ matrices over $R$. The
Lie algebra $\frak{sl}(2;\BB F_p)$ of $\SL(2;\BB F_p)$ consists of
traceless matrices: $\frak{sl}(2;\BB F_p)=\{X\in M(2;\BB
F_p)~|~X_{11}=-X_{22}\}$.  The adjoint action of $\SL(2;\BB F_p)$ on
$\frak{sl}(2;\BB F_p)$ is the action by conjugation. As a vector
space $\frak{sl}(2;\BB F_p)$ has a natural $\SL(2;\BB
F_p)$--invariant bilinear form, the Killing form $B$ defined by
$B(X,Y)=\tr(XY)$. The associated quadratic form $Q_B$ (defined by
$B(X,X)=2Q_B(X)$) is thus also invariant. Explicitly, for $X,Y\in
\frak{sl}(2;\BB F_p)$:
\[B(X,Y)=2X_{11}Y_{11}+X_{12}Y_{21}+X_{21}Y_{12}\,,\qquad
Q_B(X)=X_{11}^2+X_{12}X_{21}\,.\]

\begin{lemma}\label{le:lie}
There is a short exact sequence
\[ 1\lra \frak{sl}(2;\BB F_p)\lra \SL(2;\Z/p^2\Z)\lra \SL(2;\BB
F_p)\lra 1.\] 
The conjugation action of $\SL(2;\BB F_p)$ on $\frak{sl}(2;\BB F_p)$
induced by this sequence is the adjoint action.
\end{lemma}

\begin{proof}
  The inclusion $\Z/p\Z\to \Z/p^2\Z$ is given by $a\mapsto pa$. It
  induces an inclusion $M(2;\Z/p\Z)\to M(2;\Z/p^2\Z)$ given by
  $X\mapsto pX$.

  Reduction modulo $p$ induces the surjection $\pi\colon
  \SL(2;\Z/p^2\Z)\to \SL(2;\BB F_p)$ whose kernel is clearly
  \[\ker(\pi)=\{I+pX\in M(2;\Z/p^2\Z)~|~\det(I+pX)=1\}\,.\] Now
  $\det(I+pX)=1+p\tr(X)+p^2 \det(X)=1+p\tr(X)$ since we are in
  $\Z/p^2\Z$, so we can rewrite:
  \[\ker(\pi)=\{I+pX~|~X\in \frak{sl}(2;\BB F_p)\}\,.\]
  The equation $(I+pX)(I+pY)=I+pX+pY$ now shows that the map $X\to I+pX$
  is an isomorphism of the additive group $\frak{sl}(2;\BB F_p)$ to
  $\ker(\pi)$. The final sentence of the lemma is clear.
\end{proof}

\begin{lemma}\label{le:classes}
  The number of $\SL(2;\Z/p^2\Z)$--conjugacy classes in
  $\frak{sl}(2;\BB F_p)$ (ie, orbits of the adjoint action of
  $\SL(2;\BB F_p)$) is exactly $(p+2)$, as listed in the following
  table.  In the table $n$ represents a fixed quadratic non-residue in
  $\BB F_p$ and ``qr'' is short for quadratic residue (i.e., a
  square). Each of rows 2 and 3 represents $(p-1)/2$ conjugacy classes,
  as $Q=Q_B(X)$ runs respectively through the quadratic residues and
  non-residues in $\BB F_p-\{0\}$.
  \begin{center}
  \begin{tabular}{|c|c|c|c|}
    \hline
\emph{description}&\emph{size}&\emph{\# classes}&\emph{representative}\\
\hline\hline
trivial&$1$&$1$&$\left({0~0\atop 0~0}\right)$\\
\hline
anisotropic qr&$p(p+1)$&$(p-1)/2$&$\left({0~1\atop Q~0}\right)$\\
\hline
anisotropic non-qr&$p(p-1)$&$(p-1)/2$&$\left({0~1\atop Q~0}\right)$\\
\hline
isotropic qr&$(p^2-1)/2$&$1$&$\left({0~0\atop 1~0}\right)$\\
\hline
isotropic non-qr&$(p^2-1)/2$&$1$&$\left({0~0\atop n~0}\right)$\\
\hline
  \end{tabular}
  \end{center}
\end{lemma}

\def\Pm[#1,#2;#3,#4]{\begin{scriptsize}\begin{pmatrix}
  #1&#2\\#3&#4\end{pmatrix}\end{scriptsize}}

\begin{proof} We will prove this in several steps.
\begin{Step}
  Any $\Pm[x,y;z,-x]\in\frak{sl}(2;\BB F_p)$ is $\SL(2;\BB
  F_p)$--equivalent to a matrix of the form $\Pm[0,y';z',0]$.
\end{Step}

To see this, note first that
\[\Pm[a,b;c,d]\Pm[x,y;z,-x]\Pm[d,-b;-c,a]=\Pm[(1+2bc)x+bdz-acy,*;
*,*]\,,\]
so we want to solve the equations $ad-bc=1$ and $(1+2bc)x+bdz-acy=0$
for $a,b,c,d$.
\begin{itemize}
\item If $y\ne 0$ choose $b=0$, $a=d=1$ and solve $x-cy=0$ for $c$.
\item If $y=0$ and $z\ne 0$ choose $a=0$, $b=-c=1$ and solve for $d$.
\item If $y=z=0$ choose $2bc=-1$, $a=1$ and solve for $d$.
\end{itemize}

\begin{Step}
  If $Q=Q\Pm[x,y;z,-x]\ne 0$ then $\Pm[x,y;z,-x]$ is $\SL(2;\BB F_p)$--equivalent to
  $\Pm[0,1;Q,0]$.
\end{Step}

 We have shown we can assume $x=0$. Then
 \begin{equation}
   \label{eq:conj}
\Pm[a,b;c,d]\Pm[0,y;z,0]\Pm[d,-b;-c,a]=\Pm[bdz-acy,a^2y-b^2z;d^2z-c^2y,acy-bdz].\end{equation}
Since $Q=yz\ne 0$ we have $y,z\ne0$ so $a^2y-b^2z=1$ can be solved
for $a,b$. Then $\Pm[a,b;c,d]=\Pm[a,b;bz,ay]$ does what is required.
\begin{Step}
  Excluding the zero-element, if $Q\Pm[x,y;z,-x]= 0$ then $\Pm[x,y;z,-x]$
  is $\SL(2;\BB F_p)$--equivalent to exactly one of $\Pm[0,0;1,0]$ or $\Pm[0,0;n,0]$,
  where $n$ is a fixed quadratic non-residue.
\end{Step}

We can assume $x=0$. If $z=0$ we conjugate by an element with $a=0$ to
get $y=0$. Thus we can assume $x=y=0$ and $z\ne0$.  Now looking at
equation \eqref{eq:conj}, one sees that if $x=y=0$ then $z$
can be changed only by squares.
\begin{Step}
It remains to verify the sizes of the conjugacy classes.
\end{Step}

For each class in row 2 or 3 we must simply count the number of
elements $\Pm[x,y;z,-x]$ with $x^2+yz=Q$. Here $Q\ne 0$.  If $Q$ is a
quadratic non-residue then we must have $yz\ne 0$, so for each of $p$
choices of $x$ and each of $p-1$ choices of $y\ne 0$ we get a unique
$z$. There are therefore $p(p-1)$ elements in the class. A similar
count gives $p(p+1)$ elements if $Q$ is a residue.

If $Q=0$ it is easier to work out the isotropy group of a
representative of the class. For an element in our normal form $x=y=0$
the isotropy group consists of all $\Pm[d^{-1},0;c,d]$ with $d^2=1$.
This clearly has size $2p$ so the class has size $|\SL(2;\BB
F_p)|/2p=(p^2-1)/2$.
\end{proof}

We now investigate the $\SL(2;\BB F_p)$--classes of proper non-trivial
subgroups in $\Fr{sl}(2;\BB{F}_p)$. The group $\Fr{sl}(2;\BB{F}_p)$
itself has order $p^3$.

We first consider the subgroups of order $p$.  Using Lemma \ref{le:classes}
it is clear there are three classes. Namely
\begin{enumerate}
\item[$I$.] \emph{Isotropic lines}. Each isotropic line has $(p-1)/2$
  isotropic qr elements and $(p-1)/2$ isotropic non-qr elements. There
  are $p+1$ such lines in this class. A  representative is the
  line $\bigl\{\Pm[0,0;y,0]~|~y\in\BB F_p\bigr\}$.
\item[$R$.] \emph{Anisotropic qr lines}. Each such line has exactly two
  elements in each  anisotropic qr conjugacy class. There are
  $p(p+1)/2$ such lines in this class. A  representative is the
  line $\bigl\{\Pm[0,y;y,0]~|~y\in\BB F_p\bigr\}$.
\item[$N$.] \emph{Anisotropic non-qr lines}. Each such line has exactly two
  elements in each  anisotropic non-qr conjugacy class. There are
  $p(p-1)/2$ such lines in this class. A  representative is the
  line $\bigl\{\Pm[0,y;ny,0]~|~y\in\BB F_p\bigr\}$.
\end{enumerate}

Next, we determine the conjugacy classes of subgroups of order
$p^2$, i.e., planes. Since the Killing form $B$ is nonsingular, the
orthogonal complement of such a plane with respect to $B$ will be a
line, and vice versa, so we can classify planes up to conjugacy by the
conjugacy classes of their orthogonal complements. There are therefore
three classes of planes:

\begin{enumerate}
\item[$I^\perp$.] \emph{Orthogonal complements of isotropic lines}.
 A representative such plane is
  \[I^\perp=\bigl\{\Pm[x,0;y,-x]~|~x,y\in \BB F_p\bigr\}\,.\] The
  Killing form is degenerate on this plane, with nullspace $I$. This
  nullspace contains all isotropic elements of the plane and the
  remaining elements consist of $2p$ elements from each anisotropic qr
  conjugacy class. The plane has no anisotropic non-qr elements.
  There are $p+1$ of these planes.
\item[$R^\perp$.] \emph{Orthogonal complements of anisotropic qr
    lines}.  A representative such plane is
\[R^\perp=\bigl\{\Pm[x,y;-y,-x]~|~x,y\in \BB F_p\bigr\}\,.\] Such a plane
  has exactly $2p-2$ isotropic elements, which, together with $0$,
  form two isotropic lines (in $R^\perp$ the lines $x=y$ and $x=-y$).
  For any $Q\ne 0$ there are exactly $p-1$ elements $X\in R^\perp$
  with $Q_B(X)=Q$. Thus such a plane intersects every conjugacy
  class in $\frak{sl}(2;\BB F_p)$. There are $p(p+1)/2$ of these planes.
\item[$N^\perp$.] \emph{Orthogonal complements of anisotropic non-qr
    lines}.  A representative such plane is \[N^\perp=\bigl\{\Pm[x,y;-ny,-x]~|~x,y\in \BB F_p\bigr\}\,.\] Such a plane
  has no isotropic elements and for any $Q\ne 0$ it has $p+1$ elements with
  $Q_B(X)=Q$. There are $p(p-1)/2$ of these planes.
\end{enumerate}

We note for future reference

\begin{lemma}\label{le:elwiseconj}
Any plane of type $R^\perp$ is elementwise conjugate in
$\SL(2;\Z/p^2\Z)$ to $\frak{sl}(2;\BB F_p)$\qed
\end{lemma}

\begin{proof}[Proof of Theorem \ref{otherpair} (i)]
  We will show that the trivial subgroup is fixed point equivalent to
  any anisotropic qr line $R$. It suffices to show that
  the only finite dimensional representation of $\SL(2;\Z/p^2\Z)$
  without an $R$--fixed vector is the trivial representation. Given
  such a representation, each subgroup $H$ containing $R$ will also
  have no fixed vector. We will use this information for the subgroups
  $H$ of type $R$, $I^\perp$, $R^\perp$, $N^\perp$, and $\frak
  {sl}(2;\BB F_p)$ to show the representation must be trivial.

  To begin, the sum of the character $\chi$ of a representation over
  the non-zero elements of a line in $\Fr{sl}(2;\BB{F}_{ p})$ will
  only depend on the conjugacy class of the line, and thus give
  numbers that we shall call $X_I(\chi)$, $X_R(\chi)$, $X_N(\chi)$,
  depending on whether the line is isotropic, anisotropic qr, or
  anisotropic non-qr. Also, let $X_0(\chi)$ be the dimension of the
  representation; this is $\chi$ evaluated on the trivial element in
  $\Fr{sl}(2;\BB{F}_{p})$.  If $H$ is a subgroup of
  $\Fr{sl}(2;\BB{F}_p)$, then the sum of $\chi$ over the elements of
  $H$ gives $\abs{H}$ times the dimension of the fixed space of the
  representation restricted to $H$, hence zero under our assumption
  that $H$ has no non-trivial fixed points.  Since $H$ is a union of
  lines that are disjoint except at $0$, this then gives an equation
  of the form
  \[ X_0(\chi) + I_{H}X_I(\chi) + R_HX_R(\chi) + N_HX_N(\chi) = 0
  \,.\] Here the coefficients $I_H$, $R_H$, $N_H$ are the number of
  lines of each type in $H$. By our discussion above, these numbers for
  the subgroups of interest to us are:
\[
\begin{matrix}
  &I_H&R_H&N_H\\
H=R&0&1&0\\
H=I^\perp&1&p&0\\
H=R^\perp&2&(p-1)/2&(p-1)/2\\
H=N^\perp&0&(p+1)/2&(p+1)/2\\
H=\frak{sl}(2;\BB F_p)& p+1&(p^2+p)/2&(p^2-p)/2
\end{matrix}
\]
These five different types of subgroups containing $R$ yield five
linear equations in the four unknown quantities $X_0(\chi), X_I(\chi),
X_R(\chi), X_N(\chi)$. Since already the coefficient matrix of the
first four equations,
\[ \begin{pmatrix}
  1  &  0    &    1      &     0   \\
  1  &  1    &    p     &      0    \\
  1  &  2   &  (p-1)/2   &  (p-1)/2 \\
1 & 0  & (p+1)/2  & (p+1)/2 \\
\end{pmatrix}\,, \] has nonzero determinant (namely $-p^2$),
the equations have only the trivial solution. This implies
$X_0(\chi)=0$, proving the representation is trivial, as desired.
\end{proof}

\begin{rem}
  By computing the character table of $\SL(2;\Z/p^2\Z)$ one can show
  that there is no other fixed point equivalence in $\SL(2;\Z/p^2\Z)$
  between non-conjugate subgroups of $\frak{sl}(2;\BB F_p)$.
\end{rem}

\begin{proof}[Proof of Theorem \ref{otherpair} (ii)]
Let $\Fr{p}$ be a prime ideal of $\Cal{O}_k$ such that
$\Cal{O}_k/\Fr{p}=\BB{F}_{q}$ is a proper extension of $\BB{F}_{
  p}$ and $p$ odd; that such a prime exists follows from the Cebotarev
Density Theorem. Consider the following inclusion of short exact
sequences:
\[\xymatrix{
1\ar[r]&V_p\ar[r]\ar[d]&\PSL(2;\Z/p^2\Z)\ar[r]\ar[d]&
   \PSL(2;\BB F_p)\ar[r]\ar[d]&1\\
1\ar[r]&V_{\frak p}\ar[r]&\PSL(2;\Cal{O}_k/\Fr{p}^2)\ar[r]&
   \PSL(2;\BB F_q)\ar[r]&1\\
}\]
By Lemma \ref{le:lie} we already know the kernel $V_p$ in the first
sequence is $\frak{sl}(2;\BB F_p)$ (the transition from $\SL$ to
$\PSL$ just factors by $\{\pm I\}$ and does not affect the kernel).

Although we do not need it, we note that $V_{\frak p}=\frak{sl}(2;\BB
  F_q)$.  If $\frak p$ is principal, $\frak p=(\pi)$, say, then we
could argue as in the proof of Lemma \ref{le:lie}. In general we can
replace $k$ by its localization at $\frak p$ without changing the
second exact sequence and then $\frak p$ becomes principal, so the
argument applies.

We claim that in $\PSL(2;\Cal{O}_k/\frak p^2)$ any element of $V_p$
can be conjugated out of $V_p$. We need only show this for the
representatives of conjugacy classes given in Lemma \ref{le:classes}
and the claim is then a simple calculation using equation
\eqref{eq:conj} with $b \in \BB F_q - \BB F_p$, $a = d = 1$,
and $c = 0$.

The proof is now complete, since Lemma \ref{le:elwiseconj} gives
elementwise conjugate subgroups in $V_p$ and we have just shown they
are primitive in $\PSL(2;\Cal{O}_k/\frak p^2)$.
\end{proof}



\section{Locally symmetric manifolds and other generalities}

\label{S:Symmetric}


\subsection{$\R$--rank 1 geometries}\label{se:5.1}

 We shall denote by $\Hy^n_Y$ the $n$--dimensional hyperbolic
spaces modelled on
$Y\in\{\C,\BB{H},\BB{O}\}$ (where $n=2$ when $Y=\BB{O}$). The
methods used to produce eigenvalue,
length, and primitive length equivalent  manifolds
extend with little fuss to complex, quaternionic, and Cayley
hyperbolic manifolds. We give the version for primitive length.

\begin{thm}\label{Rrank1}
Let $\Gamma$ be a torsion-free lattice in $\Isom({\Hy}^n_Y)$. Then
there exist infinitely many pairs of finite covers of $M={\Hy}^n_Y/\Gamma$,
$\{M_j,N_j\}$ such that
\begin{itemize}
\item[(a)] $L_p(M_j)=L_p(N_j)$,
\item[(b)] $\vol(M_j)/\vol(N_j)\to \infty$.
\end{itemize}
Moreover, (a) and (b) hold for any finite volume Riemannian metric on $M$.
\end{thm}

\begin{proof}
The argument we give breaks into a few cases. First, in most cases we
have the inclusion
\[ \PO_0(B_4;\Z) < \B{G}_{Y,n}(\Z) \] where $B_4$ is the form from the
proof of Theorem \ref{T:GenPrim}, $\B{G}_{Y,n}$ is $\Q$--algebraic,
and $\B{G}_{Y,n}(\R)$ with the analytic topology is Lie isomorphic to
$\Isom(\Hy_Y^n)$. For $Y=\C$, this fails only when $n=1,2$. For
$Y=\BB{H}$, when $n\geq 3$, this is clear. The remaining cases of
$n=1,2$ follows from the exceptional isometry between $\Hy_\BB{H}^1$
and $\Hy^4$ together with the isometric inclusion of $\Hy_\BB{H}^1$
into $\Hy_\BB{H}^2$.  Finally, for $Y=\BB{O}$, this follows from the
isometric inclusion of $\Hy_\BB{H}^2$ into $\Hy_\BB{O}^2$. For all
these cases, as in the proof of Theorem \ref{T:GenPrim}, an application
of the Strong Approximation Theorem (cf \cite{We}, \cite{No}) in
combination with the Cebotarev Density Theorem provides infinitely many
primes $\Fr{p}$ such that $\Gamma$ surjects onto certain finite groups
${\bf G}(S/\Fr{p}^2)$ of Lie type which contain
$\PSL(2;\Z[i]/p^2\Z[i])$. The proof is completed just as it was in the
proof of Theorem \ref{T:GenPrim}.

It remains to deal with $Y=\C$ and $n=1,2$. The case of $n=1$ is
simply the case of hyperbolic surfaces. Case 2 cannot be handled
indirectly, and we must use primitive pairs in the finite groups
$\PU(2,1;\Cal{O}_k/\Fr{p}^2)$, where $k/\Q$ is an imaginary
quadratic extension of $\Q$ and $\Fr{p}$ is a prime ideal of
$\Cal{O}_k$. Selecting $\Fr{p}$ such that $\Cal{O}_k/\Fr{p}$ is a
quadratic extension of $\BB{F}_p$, we have the short exact sequence
\[ 1 \lra \Fr{su}(2,1;\Cal{O}_k/\Fr{p}) \lra
\SU(2,1;\Cal{O}_k/\Fr{p}^2) \lra \SU(2,1;\Cal{O}_k/\Fr{p}) \lra 1, \]
where $\Fr{su}(2,1;\Cal{O}_k/\Fr{p})$ is the Lie algebra of $\SU(2,1)$
over the field $\Cal{O}_k/\Fr{p}$. With the inclusions
\[ \Fr{sl}(2;\BB{F}_p) < \Fr{su}(2,1;\Cal{O}_k/\Fr{p}), \quad
\Omega(2,1;\BB{F}_p) < \SU(2,1;\Cal{O}_k/\Fr{p}), \]
The subgroups $\Fr{sl}(2;\BB{F}_p)$ and $R^\perp$ are elementwise conjugate
in $\SU(2,1;\Cal{O}_k/\Fr{p}^2)$ where $R^\perp$ is a
2--plane from Lemma \ref{le:elwiseconj}. It is straightforward to verify
that the pair satisfies the additional requirements needed for the
primitive case.
\end{proof}

Our methods also produce eigenvalue equivalent covers for all of these
groups as well. In addition, for sufficiently large $n$, we can
produce covers which are both primitive length and eigenvalue
equivalent; here $n\geq 5$ and $Y$ can be $\R$, $\C$, or $\BB{H}$.  To
do this, by \cite{LS} Proposition 4 Window 2 for $n\geq 5$ we can
arrange for the simple groups of orthogonal type to contain a copy of
$(\textrm{P})\SL(3;\BB{F}_p)$ which contains a group of the type given
in Theorem \ref{L:Swallow}.


\subsection{Locally symmetric manifolds}

\paragraph{Length and eigenvalue equivalent covers}

As is clear from this discussion (and the generality of the Strong
Approximation Theorem in \cite{We} and \cite{No}) our methods also
apply to lattices in every non-compact higher rank simple Lie group.
The discussion given at the end of \S~\ref{se:5.1} also applies in this setting
to arrange for the finite groups of Lie type occurring in Strong
Approximation to contain a copy of
$(\textrm{P})\SL(3;\BB{F}_p)$.

\paragraph{Primitive length equivalent covers}

 Construction of primitive length equivalent covers over a fixed locally symmetric manifold
is more subtle since in many settings the associated fundamental
group fails to have the needed condition on maximal cyclic subgroups. It
seems interesting to try to weaken the condition on maximal cyclic
subgroups to produce examples in this setting.
























\section{Final Remarks} \label{S:The End}



\subsection{Relations among length, primitive length, and eigenvalue
  equivalence}

\begin{exampleA} Let $M$ be a closed surface of genus at
least $2$ equipped with a hyperbolic metric. Let $G$ be the
alternating group $\operatorname{Alt}(4)$ and $H$ and $K$ the
elementwise conjugate pair described in \S~\ref{se:2.3}. Then, given a
surjection $\pi_1(M) \to G$, let $\gamma \in \pi_1(M)$ map to $a \in
G$ and correspond to a primitive closed geodesic in $M$ (there are
infinitely many primitive elements mapping to any element of $G$).
The non-primitive geodesic of $M$ corresponding to $\gamma^2$ has four
lifts to $M_H$, two primitive and two not, and it has three lifts to
$M_K$, all non-primitive. Of course, there might accidentally be
some unrelated primitive geodesic in $M_K$ of the right length, but
for a generic hyperbolic metric and a homomorphism to $G$ that
factors through a free group this does not happen and $M_H$ and
$M_K$ are not primitive length equivalent. Indeed, assume $\gamma$
is the shortest closed geodesic on $M$ and every other closed
geodesic has much larger length. Then if $\gamma$ maps to $a$, one
can see that $M_H$ and $M_K$ are not primitive length equivalent.
\end{exampleA}

\begin{exampleB} Eigenvalue equivalent surfaces obtained
from Theorem \ref{T:Gen2} using the trivial subspace $\set{0}$ and any
proper subspace of $\BB{F}_p^n$ generically produce examples
which are not length equivalent. In particular, eigenvalue equivalence
need not imply length or primitive length
equivalence.
\end{exampleB}

It seems plausible that length equivalent hyperbolic examples
constructed from Theorem \ref{T:Gen2} using $\BB{F}_p^n$ and any
nontrivial subspace of $\BB{F}_p^n$ will generically fail to be
eigenvalue equivalent, but this is more subtle.  Using the results of
Zelditch \cite{zelditch} it is easy to see that for a hyperbolic
manifold $M^m$ of sufficiently high dimension this approach
\emph{will} give length equivalent but not eigenvalue equivalent
examples for generic (not necessarily hyperbolic) deformations of the
metric on $M$.  Using $G=\operatorname{Alt}(4)$ this allows one to
find such examples in dimensions $m\ge 3$.

All of our examples of primitive length equivalence are also examples
of eigenvalue equivalence.

\begin{ques} Are two primitive length equivalent hyperbolic manifolds
  necessarily eigenvalue equivalent?
\end{ques}

\subsection{Complex lengths} All our results for equal length sets actually produce
manifolds which have the same \textbf{complex length sets}. Recall
that the complex length of a closed geodesic $\gamma$ in a Riemannian
$m$--manifold is a pair $(\ell(\gamma),V)$ where $\ell(g)$ is the
length of $\gamma$ and $V\in\Ort(m-1)$ is determined by the holonomy
of $\gamma$. The complex length spectrum is the collection of such
complex lengths with multiplicities, and the complex length set forgets
multiplicities as before. The point is that Theorem \ref{T:Swallow}
gives manifolds with the same complex length sets, just as Sunada's
theorem gives equal complex length spectra. See \cite{Sa} for more on
the complex length spectrum.

\subsection{Commensurability} The known methods of producing
isospectral or iso-length spectral hyperbolic manifolds result in
commensurable manifolds and it is an open question as to whether this
is always the case.  By construction, the eigenvalue and (primitive)
length equivalent hyperbolic manifolds constructed here are also
commensurable.

\begin{ques}\label{ques}
  Let $M_1$ and $M_2$ be eigenvalue (resp.\ length or primitive
  length) equivalent closed hyperbolic manifolds.  Are they
  commensurable?
\end{ques}

There has been some recent activity on this question. It is shown
that Question \ref{ques} has an affirmative answer in the length
equivalent setting if the manifolds $M_1$ and $M_2$ are arithmetic
hyperbolic 3--manifolds (\cite{CHLR}), or if the manifolds are even
dimensional arithmetic hyperbolic manifolds (\cite{PrRap}). Indeed,
the results of \cite{PrRap} apply to more general locally symmetric
spaces. In contrast, \cite{PrRap} also exhibts arbitrarily large collections
of incommensurable hyperbolic 5--manifolds which are length
commensurable. The commensurability classes of these manifolds seem
to be the best candidates for producing a negative answer to Question
\ref{ques}.

\subsection{Infinite sets of examples} Our constructions show that there can be no uniform
bound on the number of pairwise eigenvalue (resp.\ length or primitive
length) equivalent, non-isometric manifolds. Thus a natural question
is.

\begin{ques}
  Are there infinite sets of pairwise eigenvalue (resp.\ length or
  primitive length) equivalent, closed hyperbolic $m$--manifolds?
\end{ques}

In the context of length equivalence a positive answer would follow if
one can find infinitely many mutually elementwise conjugate subgroups
of finite index in a finitely generated free group. C. Praeger pointed
out to us that a slightly stronger version of this question is listed
as an open problem (Problem 11.71) in the Kourovka Notebook
\cite{kourovka}.  It was asked there in the parallel context of
Kronecker equivalence of number fields. It seems likely that the
answer to this question is ``no'', but the limited partial answers
that are known involved considerable effort, see \cite{praeger}.

\subsection{Can one hear the size of a drum?}

Mark Kac's famous paper ``Can one hear the shape of a drum''
\cite{kac} is quoted in many papers on isospectrality. Of course,
the ``drums'' of his title were not closed hyperbolic manifolds, but
rather flat plane domains.  The first pair of different ``drums''
with the same sound (i.e., non-isometric isospectral plane domains)
was found in the 1990's by Gordon, Webb, and Wolpert \cite{drums}.

However, one might question whether the sounds of their drums $D_1$
and $D_2$ are really indistinguishable. They comment: ``... to
produce the same sound (i.e., the same frequencies with the same
amplitudes) as would result from striking $D_1$ at a given point
with a given (unit) intensity ... one must strike $D_2$
simultaneously at seven points with appropriate intensities''.  A
more obvious example of this issue is a pretty example of S. Chapman
\cite{chapman}. Chapman reinterprets earlier discussion of the
Gordon--Webb--Wolpert examples in terms of paper folding and
cutting, as is familiar from making paper dolls. Of course, by
cutting too much one can create disconnected objects, and by this
means Chapman derives from the Gordon--Webb--Wolpert example the
following simple example: $D_1$ is the disjoint union of a unit
square and an isosceles right triangle with legs of length $2$, and
$D_2$ is the disjoint union of a $1\times 2$ rectangle and an
isosceles right triangle with legs of length $\sqrt2$. This pair of
domains is isospectral, but one can ask to what extent they really
sound the same.

A more honest example of equal sound might be the following:
purchase three identical drums and let $D_1$ consists of one of them
and $D_2$ consist of the disjoint union of the other two. It would
be hard to distinguish $D_1$ from $D_2$ on hearing a drummer strike
either one once.  This example suggests that eigenvalue equivalence
may have as much right as isospectrality to be interpreted as ``same
sound.''

In his paper Kac gave a proof that drums that sound the same have
equal area, but this was based on isospectrality. Revisiting this in
the context of eigenvalue equivalence we ask:

\begin{ques}
  Do there exist connected eigenvalue equivalent plane domains of
  unequal area?
\end{ques}





\noindent Department of Mathematics\\
University of Illinois at Urbana-Champaign\\
Urbana, IL 61801\\
email: {\tt clein@math.uiuc.edu}\\

\noindent Department of Mathematics\\
California Institute of Technology\\
Pasadena, CA 91125\\
email: {\tt dmcreyn@caltech.edu}\\

\noindent Department of Mathematics\\
 Barnard College, Columbia University \\
 New York, NY 10027 \\
email: {\tt neumann@math.columbia.edu} \\

\noindent Department of Mathematics\\
University of Texas\\
Austin, TX 78712\\
email: {\tt areid@math.utexas.edu}



\end{document}